\documentclass[twoside,notitlepage,11pt]{amsart}

\usepackage{amssymb}
\usepackage[leqno]{amsmath}
\usepackage{amsfonts}
\usepackage{amsopn}
\usepackage{amstext}
\usepackage{amsthm}

\providecommand{\cal}{\mathcal}
\renewcommand{\Bbb}{\mathbb}
\renewcommand{\frak}{\mathfrak}
\newenvironment{pf}{\begin{proof}}{\end{proof}}

\newcommand{\Aaa}{{\cal{A}}}
\newcommand{\Bee}{{\cal{B}}}
\newcommand{\Cee}{{\cal{C}}}

\newcommand{\Kay}{{\cal{K}}}

\newcommand{\Pee}{{\cal{P}}}

\newcommand{\Err}{{\Bbb{R}}}

\newcommand{\al}{\alpha}

\newcommand{\sig}{\sigma}

\renewcommand{\phi}{\varphi}
\renewcommand{\rho}{\varrho}

\newcommand{\rest}{\restriction}

\newcommand{\ntr}{n\in\omega}

\newcommand{\loe}{\leqslant}
\newcommand{\goe}{\geqslant}

\newcommand{\subs}{\subseteq}
\newcommand{\sups}{\supseteq}
\newcommand{\nnempty}{\ne\emptyset}

\newcommand{\cl}{\operatorname{cl}}
\newcommand{\Int}{\operatorname{int}}

\newcommand{\id}{\operatorname{id}}

\newcommand{\poset}{{\Bbb{P}}}

\newcommand{\meet}{\cdot}

\newcommand{\join}{+}

\newtheorem{tw}{Theorem}[section]
\newtheorem{wn}[tw]{Corollary}
\newtheorem{lm}[tw]{Lemma}
\newtheorem{prop}[tw]{Proposition}
\newtheorem{claim}[tw]{Claim}
\theoremstyle{definition}
\newtheorem{df}[tw]{Definition}
\newtheorem{ex}[tw]{Example}
\newtheorem{exs}[tw]{Examples}
\theoremstyle{remark}
\newtheorem{uwgi}[tw]{Remark}

\newcommand{\setof}[2]{\{#1\colon #2\}}
\newcommand{\bigsetof}[2]{\Bigl\{#1\colon #2\Bigr\}}
\newcommand{\seq}[1]{\langle #1 \rangle}

\newcommand{\sn}[1]{\{#1\}} 
\newcommand{\pair}[2]{{\langle #1, #2 \rangle}} 
\newcommand{\map}[3]{#1\colon #2 \to #3} 
\newcommand{\img}[2]{#1[#2]} 
\newcommand{\inv}[2]{{#1}^{-1}[#2]} 

\newcommand{\power}[1]{\Pee(#1)}
\newcommand{\dpower}[2]{[#1]^{#2}}

\newcommand{\fin}[1]{[#1]^{<\omega}}

\newcommand{\En}{\mathcal N}

\newcommand{\tight}{\operatorname{t}}

\newcommand{\cont}{\ensuremath{2^{\omega}}}

\newcommand{\ro}{\operatorname{RO}}

\newcommand{\elsig}{\operatorname{L\Sigma}}
\newcommand{\Elsig}{\operatorname{L\Sigma}}

\newcommand{\elsigO}{\ensuremath{\elsig(\loe\omega)}} 
\newcommand{\elsigfin}{\ensuremath{\elsig(<\omega)}} 
\newcommand{\kelsig}{\operatorname{KL\Sigma}}
\newcommand{\kelsigO}{\ensuremath{\kelsig(\loe\omega)}}

\newcommand{\MA}{\operatorname{MA}}
\newcommand{\MAone}{\ensuremath{\MA_{\omega_1}}}
\newcommand{\lat}{\ensuremath{\mathbb L}}
\newcommand{\ppi}{\ensuremath{\pi}}
\newcommand{\closed}{\operatorname{Closed}}
\newcommand{\ult}{\operatorname{Ult}}

\providecommand{\K}{{\mathbb K}}

\newtheorem{question}[tw]{Question}

\renewcommand\:{\nobreak \hskip .1111em\mathpunct {}\nonscript \mkern
-\thinmuskip {:}\hskip .3333emplus.0555em\relax}
\newcommand\Cal{\mathcal}
\begin{document}

\title{On some classes of Lindel\"of $\Sigma$-spaces}
\author{Wies{\l}aw Kubi\'s}
\address{Institut Matematyki, Uniwersitet \'Sl\c aski, ul. Bankowa 14,
40-007 Katowice, Polska}
\email{kubis@ux2.math.us.edu.pl}
\thanks{The first author's research was supported by a postdoctoral fellowship at York
University, Toronto, Canada (2003/2004). The first author would also like to
thank the Fields Institute for their great hospitality during the work on
this article.}
\author{Oleg Okunev}
\address{Facultad de Ciencias F\'\i sico-Matem\'aticas, Benem\'erita
Universidad Aut\'onoma de Puebla \\ av. San Claudio y Rio Verde s/n
col. San Manuel, Ciudad Universitaria\\ CP 72570 Puebla, Puebla, M\'exico}
\email{oleg@servidor.unam.mx}
\thanks{The second author acknowledges support from PROMEP, 103.5/04/2539}
\author{Paul J.\ Szeptycki}
\address{Mathematics, Atkinson Faculty, York University, Toronto, ON Canada M3J 1P3}
\email{szeptyck@yorku.ca}
\thanks{The third author acknowledges support from NSERC grant 238944}



\begin{abstract} We consider special subclasses of the class of Lindel\"of $\Sigma$-spaces obtained by imposing restrictions on the weight of the elements of compact covers that admit countable
networks: A space $X$ is in the class $L\Sigma(\loe\kappa)$ if it admits a cover by compact subspaces of weight
$\kappa$ and a countable network for the cover. We restrict our attention to $\kappa\leq\omega$. 
In the case $\kappa=\omega$, the class includes the class of metrizably fibered spaces considered by Tkachuk, and the $P$-approximable spaces considered by Tka\v cenko. The case $\kappa=1$ corresponds to the spaces of countable network weight, but even the case $\kappa=2$ gives rise to a nontrivial class of spaces. The relation of known classes of compact spaces to these classes is considered. It is shown that not every Corson compact of weight $\aleph_1$ is in the class $L\Sigma(\leq \omega)$, answering a question of Tkachuk. As well, we study whether certain compact spaces in $L\Sigma(\leq\omega)$ have dense metrizable subspaces, partially answering a question of Tka\v cenko. Other interesting results and examples are obtained, and we conclude the paper with a number of open questions.

\end{abstract}

\maketitle

\section{Preliminaries}

All spaces we consider are Tychonoff (that is, completely regular
Hausdorff), unless otherwise indicated. We use terminology and
notation as in \cite{engelking}, with the exception that the tightness
of a space $X$ is denoted as $t(X)$. 

Given a locally compact space $X$, we denote by $\al X$ its one-point
compactification, the new point will be usually denoted by $\infty$
(unless $\infty\in X$). The one-point compactification of a discrete
space of size $\kappa$ will be denoted by $A(\kappa)$.

A topological space $X$ is a {\em Lindel\"of $\Sigma$-space}, if it
is a $\Sigma$-space in the sense of Nagami \cite{nag} and is Lindel\"of.
Lindel\"of $\Sigma$-spaces are also known as {\em K-countably
determined spaces} \cite{rj}.

A {\it multivalued mapping} from $X$ to $Y$ is a mapping that assigns
to every point of $X$ a subset of $Y$ (not necessarily nonempty).
For a multivalued mapping $p\:X\to Y$ and a set $A\subset X$, the {\em
image} of $A$ under $p$ is the set 
$$
p(A)=\bigcup\{p(x): x\in A\};
$$
a mapping $p:X\to Y$ is {\em onto} if $p(X)=Y$. If $p\:X\to Y$ and
$q\:Y\to Z$ are multivalued mappings, then the {\em composition} of
$q$ and $p$ is the multivalued mapping $q\circ p\:X\to Z$ such that
$(q\circ p)(x)=q(p(x))$ for all $x\in X$. 

We always use adjectives such as ``compact-valued'', ``finite-valued''
etc. for multivalued mappings; the word ``mapping'' (or ``function'')
without such an adjective will always mean a usual single-valued mapping
(which we naturally identify with the corresponding singleton-valued
mappings).

A multivalued mapping $p\:X\to Y$ is called {\em upper semicontinuous}
(or {\it usc}) if for every open set $V$ in $Y$, the set $\{\,x\in X:
p(x)\subset V\,\}$ is open in $X$. It is easy to see that
continuous functions and inverses of perfect mappings are
compact-valued usc. If $F$ is a closed subspace of $X$, then the
mapping $p_F\: X\to F$ defined by
$$
p_F(x)=\begin{cases} \{x\}\quad &\text{if $x\in F$}\\
\,\,\,\emptyset&\text{if $x\notin F$}\end{cases}
$$
(the inverse of the embedding $i_F\:F\hookrightarrow X$) is usc. A
straightforward verification shows that a composition of
compact-valued usc mappings is compact-valued and usc; a standard
argument using the closedness of the graph of $p$ in $X\times \beta Y$
proves the following:

\begin{prop} Let $p\:X\to Y$ be a multivalued mapping. Then the
following conditions are equivalent:

\begin{enumerate}\label{p0}
\item[$(a)$] $p$ is compact-valued usc.
\item[$(b)$] $p$ is a composition of the inverse of a perfect mapping onto a
closed subspace of $X$ and a continuous function.
\item[$(c)$] There are a compact space $K$, a closed subspace $F$ of\/
$X\times K$ and a continuous function $f\:F\to Y$ such that $p=f\circ
i_F^{-1}\circ\pi_X^{-1}$, where $\pi_X\:X\times K\to X$ is the
projection and $i_F\:F\to X\times K$ is the embedding.
\end{enumerate}
\end{prop} 

A family of sets $\Cal N$ is called a {\em network with respect to
a cover $\Cal C$} of a space $X$ if for every set $C\in C$ and every
neighborhood $U$ of $C$ there is an element $N$ of $\Cal N$ such that
$C\subset N\subset U$ \cite{nag}.

Note that if $\Cal C$ is a compact cover of $X$ (that is, all elements
of $\Cal C$ are compact), and $\Cal N$ is a network with respect to
$\Cal C$, then the family of the closures of the elements of $\Cal N$
is also a network with respect to $\Cal C$. 
 
The next proposition sums up several well-known characterizations of
Lindel\"of $\Sigma$-spaces \cite{arh-ci}, \cite{rj}.

\begin{prop}\label{p1} Let $X$ be a space. The following
conditions are equivalent:

\begin{enumerate}
	\item[$(a)$] $X$ is a Lindel\"of\/ $\Sigma$-space.
	\item[$(b)$] There are a compact cover $\Cal C$ of $X$ and a
countable network $\Cal N$ with respect to $\Cal C$.
        \item[$(c)$] There are a second-countable space $M$ and a
compact-valued usc mapping $p\:M\to X$ such that $p(M)=X$.
	\item[$(d)$] There are a second-countable space $M$, a space
$L$ and mappings $g\:L\to M$ and $f\:L\to X$ such that $g$ is perfect
and $f$ is continuous onto.
        \item[$(e)$] There are a second-countable space $M$, a compact
space $K$, a closed subspace $F$ of\/ $M\times K$ and a continuous
mapping $f\:F\to X$ such that $f(F)=X$.
\end{enumerate}
\end{prop}

The equivalence of $(a)$ and $(b)$ is immediate from the definition
(see \cite{nag}). The equivalence of $(c)$, $(d)$ and $(e)$ follows
from Proposition \ref{p0}. If $p\:M\to X$ is a compact-valued usc
mapping onto $X$, and $\Cal B$ is a countable base for $M$, then 
$\{\,p(B):B\in \Cal B\,\}$ is a countable network with respect to the
compact cover $\{\,p(m):m\in M\,\}$ of the space $X$, so $(c)$ implies
$(b)$. To verify that $(b)$ implies $(c)$, equip $\Cal N$ with the
discrete topology, and let $M$ be the subspace of $\Cal N^\omega$
(equipped with the product topology) consisting of all functions
$m:\omega\to\Cal N$ with the property that
$\{\,m(i):i\in\omega\,\}=\{\,N\in \Cal N:C\subset N\,\}$ for some
$C\in\Cal C$. Then $M$ is a second-countable space; the mapping
$p\:M\to X$ defined by the rule 
$$
p(m)=\bigcap\{\,m(i):i\in\omega\,\}
$$
is compact-valued, usc, and onto $X$.

\medskip
Note that the cardinality of the cover $\Cal C$ as in $(b)$ cannot
exceed $2^\omega$.

\bigskip

Of course, all compact spaces and all spaces with a countable network
are Lindel\"of $\Sigma$-spaces. From Proposition \ref{p1} it follows easily
that the class of Lindel\"of $\Sigma$-spaces is invariant with
respect to images under compact-valued usc mappings (in particular,
continuous images, closed subspaces and perfect preimages), 
countable products and countable unions.

\section{The classes $\elsig(\loe\kappa)$ and $\kelsig(\loe\kappa)$}

\medskip
In this article we consider subclasses of the class of all Lindel\"of
$\Sigma$-spaces obtained by requiring that the elements of the compact
cover $\Cal C$ as in Proposition~\ref{p1}(b) have a given
property. This leads to the following definition.

\begin{df} Let $\Kay$ be a class of compact spaces. Define
$\elsig(\Kay)$ as the class of all spaces such that there are a
second-countable space $M$ and a compact-valued usc mapping $p\:M\to X$
such that $p(M)=X$ and $p(m)\in\Kay$ for all $m\in M$.

We also define the class $\kelsig(\Kay)$ as the class of all spaces
such that there are a {\em compact\/} second-countable space $M$ and a
compact-valued usc mapping $p\:M\to X$ such that $p(M)=X$ and
$p(m)\in\Kay$ for all $m\in M$.
\end{df}

Clearly, always $\kelsig(\Kay)\subs\elsig(\Kay)$ and all spaces in
$\kelsig(\Kay)$ are compact.

\medskip

An argument similar to the proof of Proposition~\ref{p1} gives the
following:

\begin{prop}\label{p2-0} Let $X$ be a space and $\Kay$ a class of
compact spaces. Then the following conditions are equivalent:

\begin{enumerate}
	\item[$(a)$] $X\in \elsig(\Kay)$.
	\item[$(b)$] There are a compact cover $\Cal C$ of $X$ such that
$\Cal C\subset \Kay$ and a countable network $\Cal N$ with
respect to $\Cal C$.
\end{enumerate}

If the class $\Kay$ is closed with respect to continuous images and
closed subspaces, then these conditions are also equivalent to
\begin{enumerate}
	\item[$(c)$] There are a second-countable space $M$, a space
$L$ and mappings $g\:L\to M$ and $f\:L\to X$ such that $g$ is perfect,
$f$ is continuous, and $g^{-1}(m)\in\Kay$ for all $m\in M$.
        \item[$(d)$] There are a second-countable space $M$, a compact
space $K$, a closed subspace $F$ of\/ $M\times K$ and a continuous
mapping $f\:F\to X$ such that $f(F)=X$ and
$F\cap \pi_M^{-1}(m)\in \Kay$ for all $m\in M$, where $\pi_M\:M\times
K\to M$ is the projection.
\end{enumerate}
\end{prop}

Similarly,

\begin{prop}\label{p2-1} Let $X$ be a space and $\Kay$ a class of
compact spaces is invariant under continuous images and
closed subspaces. Then the following conditions are equivalent:

\begin{enumerate}
	\item[$(a)$] $X\in \kelsig(\Kay)$.
	\item[$(b)$] There are a compact second-countable space $M$, a
space $L$ and mappings $g\:L\to M$ and $f\:L\to X$ such that $g$ is
perfect, $f$ is continuous, and $g^{-1}(m)\in\Kay$ for all $m\in M$.
        \item[$(c)$] There are a compact second-countable space $M$, a
compact space $K$, a closed subspace $F$ of\/ $M\times K$ and a
continuous mapping $f\:F\to X$ such that $f(F)=X$ and
$F\cap \pi_M^{-1}(m)\in \Kay$ for all $m\in M$.
\end{enumerate}
\end{prop}

\medskip
Thus, a compact space $X$ is in $\elsig(\Kay)$ if and only if $X$ has
a countable closed cover $\Cal N$ such that for every $x\in X$ the set
$\bigcap\{\,N\in\Cal N:x\in N\,\}$ belongs to $\Kay$. A (not necessarily
compact) space $X$ satisfying this condition is called {\em
$\Kay$-approximable} in \cite{tkacenko} and {\em weakly
$\Kay$-fibered} in \cite{tkachuk} (in fact, \cite{tkachuk} deals only
with the class of {\it weakly metrizably fibered\/} spaces which is
the class of weakly $\Kay$-fibered spaces with $\Kay$ the class of
all metrizable compacta). Note that a countably
compact space which is $\Kay$-approximable is in $L\Sigma(\Kay)$.

\medskip
If $\kappa$ is a cardinal, finite or infinite, we denote by
$\elsig(\loe\kappa)$ and  $\kelsig(\loe\kappa)$ the classes
$\elsig(\Kay)$ and $\kelsig(\Kay)$ where $\Kay$ is the class
of all compact spaces of weight $\loe\kappa$. Similarly,
$\elsig(<\kappa)$ and  $\kelsig(<\kappa)$ are the classes 
$\elsig(\Kay)$ and $\kelsig(\Kay)$ where $\Kay$ is the class
of all compact spaces of weight $<\kappa$; let
$\elsig(\kappa)=\elsig(\loe\kappa)\setminus\elsig(<\kappa)$.
Since the cardinality of a compact cover with respect to which there
is a countable network does not exceed $2^\omega$, all spaces in 
$\elsig(\loe\kappa)$ have cardinality at most $2^{\kappa+\omega}$, and if $\kappa\goe 2^\omega$, then the class $L\Sigma(\loe\kappa)$ coincides with
the class of all Lindel\"of $\Sigma$-spaces of network weight $\loe\kappa$.

When $\kappa$ is a finite cardinal, ``$\loe\kappa$" means ``at most
$\kappa$-element sets". Thus, $X\in\elsig(n)$, $n\in\omega$, if $X$
has a cover $\Cee$ consisting of at most $n$-element sets which has a
countable network in $X$, but $X$ does not have such a cover
consisting of at most $(n-1)$-element sets. Obviously, $\elsig(\loe
1)$ is the class of all spaces of countable network weight, and
$\kelsig(\loe 1)$ is the class of all metrizable compacta.

>From Propositions \ref{p2-0} and \ref{p2-1} readily follows

\begin{prop}\label{ls-invariance} Let $\kappa$ be a cardinal.
Then the classes $\elsig(\loe\kappa)$, $\elsig(<\kappa)$, 
$\kelsig(\loe\kappa)$ and $\kelsig(<\kappa)$ are invariant with
respect to closed subspaces, continuous images and finite unions. The
classes $\elsig(\loe\kappa)$ and $\elsig(<\kappa)$ are invariant with
respect to countable unions.
\end{prop} 

Since the product of a family of compact-valued usc mappings is
compact-valued and usc, we have

\begin{prop}\label{ls-products} If $X\in \elsig(\loe\kappa)$ and
$Y\in\elsig(\loe\lambda)$, then $X\times Y\in \elsig(\loe\lambda\cdot\kappa)$.
\end{prop}

and

\begin{prop}\label{ls-products-general} If $\{X_n:n\in\omega\,\}$
is a countable family of spaces, and $X_n\in\elsig(\loe\kappa_n)$,
$n\in \omega$, then $\prod\{\,X_n:n\in \omega\,\}\in \elsig(\loe\kappa)$
where $\kappa=|A|\cdot\sup\{\kappa_n:n\in \omega\,\}$.
\end{prop}

\medskip

\begin{exs}\label{exs-basic} 1. The double arrow space is in $\kelsig(2)$. Indeed, it
admits a 2-to-1 perfect mapping onto the closed interval, and
therefore is in $\kelsig(\loe 2)$.  It is not in $\kelsig(1)$, because
it has no countable network.

\medskip
2. Let $T$ be the unit circle, and $AD(T)$ its Alexandroff
duplicate (see e.g. \cite{engelking}). Then $AD(T)$ is in
$\kelsig(2)$, because $AD(T)$ is not metrizable and admits a perfect
2-to-1 mapping onto T.

\medskip
3. The space $A(2^\omega)$ is a non-metrizable continuous image of
$AD(T)$, and hence $A(2^\omega)\in\elsig(2)$. Therefore, 
$A(\kappa)\in\elsig(\loe\omega)$ iff $A(\kappa)\in\elsig(\loe 2)$
iff $\kappa\loe 2^\omega$.
\end{exs}

A compact space $X$ is called {\em metrizably fibered} \cite{tkachuk}
if $X$ admits a continuous mapping with metrizable fibers onto a
metrizable space. Clearly, all metrizably fibered spaces are in
$\kelsigO$, and by Proposition \ref{p2-1} every space in
$\kelsig(\loe\omega)$ is a continuous image of a metrizably fibered
compact space. Since all metrizably fibered compact spaces are
first-countable, $A(\omega_1)$ is in $\kelsig(2)$, but is not metrizably
fibered. Note also that every space in $\elsig(\loe\omega)$ is weakly
metrizably fibered, and that a compact space is weakly metrizably
fibered if and only if it is in $\elsig(\loe\omega)$ (the family of
all finite intersections of members of the family $\Cal N$ from the
definition of weakly $\Kay$-metrizable spaces cited above is a network
with respect to the cover of $X$ formed by the sets 
$C_x=\bigcap\{\,N\in\Cal N:x\in N\,\}$, $x\in X$). It is shown in
\cite{GSz} that every such compact space is sequential.

Since every metrizably fibered compact space is first-countable, it
follows that every space in $\kelsig(\loe\omega)$ is Fr\'echet.
 
\medskip
\begin{ex} (Example 2.13 in \cite{tkachuk}). Let $K$ be
the one-point compactification of a Mr\'owka space. Then $K$ is a
countable union of subspaces in $\elsig(\loe2)$, (countably many
singletons and the one-point compactification of a discrete space of
cardinality $\loe2^\omega$) and hence is itself in $\elsig(\loe2)$.
Since $K$ is not Fr\'echet, it is not in $\kelsig(\loe\omega)$.
\end{ex}

As we mentioned above, every compact space in $\elsig(\loe\omega)$
is sequential \cite{GSz}, and therefore has countable tightness
\cite{tkachuk}. Of course, every space in $\elsig(\loe 1)$ has
countable tightness. The next example shows that not all spaces in
$\elsig(<\omega)$ have countable tightness.

\begin{ex} Let $X$ be the subspace of $2^{\omega_1}$ which is the union
of the set $S$ of all points that have finitely many coordinates equal
to 1 and the singleton \{{\bf 1}\} where {\bf 1} is the point whose
all coordinates are equal to 1. It is easy to see that the tightness
of $X$ at the point {\bf 1} is uncountable. From Lemma~2.9 in
\cite{ok-tka} it follows that $S$ is a countable union of continuous
images of spaces of the form $A(\omega_1)^n\times 2^n$,
$n\in\omega$, and hence is in the class $\elsig(<\omega)$. Thus, $X$
is a $\sigma$-compact space in $\elsig(<\omega)$ of uncountable
tightness. 
\end{ex}

\begin{question} Does the class $\elsig(2)$ contain a space of
uncountable tightness? Does any of the classes $\elsig(n)$,
$n\in\omega$, contain a space of uncountable tightness?
\end{question}

\medskip
Recall that a {\em free sequence of length $\kappa$} in a topological
space $X$ is a function $f\:\kappa\to X$ such that for every
$\alpha<\kappa$, the sets $\{\,f(\beta):\beta<\alpha\,\}$ and 
$\{\,f(\beta):\alpha\loe\beta<\kappa\,\}$ have disjoint closures. 
If $X$ is a compact space, then the tightness of $X$ is equal to the
supremum of the lengths of free sequences in $X$ \cite{ava-freeseq}.

\begin{tw}\label{ciasnota} Assume that $\kappa$ is an uncountable
regular cardinal, and $X\in\elsig(<\kappa)$. Then every free sequence
in $X$ has length $<\kappa$.
\end{tw}

\begin{pf} Fix a compact cover $\Cee$ with a countable network $\En$
in $X$ so that every element of $\Cee$ has weight $<\kappa$.

Suppose $\map f{\kappa}X$ is a free sequence in $X$. For every
$\alpha<\beta\loe\kappa$ put $F(\al,\beta)=\cl\img
f{[\al,\beta)}$. Then by the definition of a free sequence 
$F(0,\al)\cap F(\al,\kappa)=\emptyset$. 

Since $|\En|=\omega$, there exists $\delta<\kappa$ such that for
every $N\in\En$ either $\sup\inv fN<\delta$ or $\inv fN$ is unbounded
in $\kappa$. Fix $C\in\Cee$ so that $f(\delta)\in C$. We claim that
there exists $\al_0\in[\delta,\kappa)$ such that $C\cap
F(\al_0,\beta)=\emptyset$ for every $\beta>\al_0$. Indeed, otherwise
for every $\al>\delta$ pick $p_\al\in C$ such that $p_\al\in
F(\al,\rho(\al))$ for some $\rho(\al)>\al$. Choose unbounded
$S\subs\kappa\setminus \delta$ such that $\rho(\al)<\al'$ whenever
$\al,\al'\in S$ and $\al<\al'$. Then $\setof{p_\al}{\al\in S}$ is a
free sequence of length $\kappa$ in $C$, which contradicts the assumption
that $w(C)<\kappa$.

Fix $\alpha_0>\delta$ so that $C\cap F(\al_0,\beta)=\emptyset$ for
every $\beta>\al_0$. Find $\beta_0>\al_0$ such that $N\cap \inv
f{[\al_0,\beta_0)}\nnempty$ whenever $\inv fN$ is unbounded in
$\kappa$ and $N\in\En$. Then $X\setminus F(\al_0,\beta_0)$ is a
neighborhood of $C$. Thus, there exists $N\in\En$ such that $C\subs N$
and $N\cap F(\al_0,\beta_0)=\emptyset$. Then $\inv fN$ is bounded in
$\kappa$ and therefore $\sup \inv fN<\delta<\al_0$. On the other hand
$f(\delta)\in N$, a contradiction.
\end{pf}

\begin{wn}\label{tightness} Assume $X\in\elsig(\loe\kappa)$ is compact and
$\kappa\goe\omega$. Then $\tight(X)\loe\kappa$. \end{wn}

\begin{uwgi} The same argument proves the following statement: 
{\sl Assume that $\kappa$ is an uncountable regular cardinal, and
$X\in\elsig(t<\kappa)$. Then every free sequence in $X$ has length
$<\kappa$}; here ``$t<\kappa$'' is the class of all compact spaces of
tightness $<\kappa$. A similar statement for countably compact spaces
was proved by Tka\v cenko \cite[Assertion 2.2]{tkacenko}. A special case
(for compact spaces) was proved by Tkachuk [13, Thm. 2.11]
\end{uwgi}

\medskip
Gerlits and Szentmikl\'ossy proved in \cite{GSz} that the Helly space
belongs to $\elsigO$ (in fact they showed that the Helly space can be
mapped onto a metric space by a map with metrizable fibers -- this is
a stronger property than being in $\elsigO$). Todor\v cevi\'c proved
in \cite{Todorcevic1999} a dichotomy for Rosenthal compact spaces,
where one of the assertions is ``being a two-to-one preimage of a
compact metric space". In view of the next result, it is natural to
ask whether all Rosenthal compacta belong to $\elsigO$.

\begin{prop} For every Polish space $X$, the space 
$$
DC_\omega(X)=\setof{f\in \Bbb R^X}{f\text{ has only countably
many points of discontinuity }}
$$
is metrizably-approximable. In particular, every compact subspace of
$DC_\omega(X)$ belongs to $\elsigO$. \end{prop}

\begin{pf} Fix a countable base $\Bee$ in $X$. Define
$$
N_{u,J}=\setof{f\in DC_\omega(X)}{\img f{\cl u}\subs J}
$$
and let $\En=\setof{N_{u,J}}{u\in\Bee,\;J\text{ is a closed rational
interval}}$. Then $\En$ is a countable family of closed subsets of
$DC_\omega(X)$.

Let $C_f=\bigcap\setof{N\in\En}{f\in N}$. We claim that $C_f$ is
metrizable for every  $f\in DC_\omega(X)$. 

Fix $f\in DC_\omega(X)$ and $g\in C_f$. We have $g(x)=f(x)$
whenever $f$ is continuous at $x$. Indeed, if $f(x)\ne g(x)$, then we
can find $u\in \Bee$ and a closed rational interval $J$ such that
$x\in u$, $\img f{\cl u}\subs J$ and $g(x)\notin J$. Then $f\in
N_{u,J}\in\En$ and $g\notin N_{u,J}$, a contradiction.

Let $A$ denote the set of all points of discontinuity of $f$. We have
proved that
$$
C_f\subs\setof{h\in\Err^X}{f\rest(X\setminus A)=h\rest(X\setminus
A)}.
$$

The set on the right-hand side is homeomorphic to $\Err^A$ with
the product topology. It follows that $C_f$ is metrizable, because
$|A|\loe\omega$.
\end{pf}

\section{Small Corson compacta need not be in $\elsigO$}

Tkachuk proved in \cite{tkachuk} that every Eberlein compact space of
weight at most continuum (equivalently, of size at most continuum)
belongs to $\elsigO$. This is not true for Corson compacta. A
consistent counterexample is given in \cite{tkachuk}. We show that a
certain Corson compact space constructed (in ZFC) by Todor\v cevi\'c
in \cite{Todorcevic1981} (see also \cite[p. 287]{Todorcevic1984}) does
not belong to $\elsigO$. This answers Tkachuk's question from
\cite{tkachuk}.

We recall the construction of Todor\v cevi\'c's Corson compact
space. Fix a stationary co-sta\-tio\-na\-ry set $A\subs\omega_1$ and
let $T(A)$ be the collection of all subsets of $A$ which are closed in
$\omega_1$. Then $\pair{T(A)}{\subs}$ is a tree of height $\omega_1$
whose all branches are countable. Let $Y(A)$ denote the collection of
all initial branches of the tree $T(A)$ (an {\em initial branch} is a
linearly ordered subset which is also closed downwards). Then $Y(A)$
is a compact subspace of the Cantor cube $\power{T(A)}$. Since all
branches of $T(A)$ are countable, $Y(A)$ is Corson compact.

Note that $|T(A)|\loe{\omega_1}^{\omega}=2^{\omega}$ and therefore
$Y(A)$ is a ``small" Corson compact space. 

\begin{prop} For every stationary co-stationary set $A\subs\omega_1$,
$Y(A)\notin\elsigO$.
\end{prop}

\begin{pf}
Suppose $Y(A)\in\elsigO$ and let $\Cee$ be a cover of $Y(A)$
consisting of metric compacta and $\En$ a countable network
with respect to $\Cee$. It has been proved in \cite[p. 287]{Todorcevic1984}
that $T(A)$ is a Baire partial order, which means that, as a forcing
notion, it does not add new countable sequences. 

Assume now that we are working in a countable transitive ZFC model $M$
and let $G$ be a $T(A)$-generic filter over $M$. In $M[G]$, the space
$Y(A)$ (with the topology generated by open sets from the ground
model) is still in $\elsigO$, because $\Cee$ still consists of metric
compacta (by the fact that there are no new sequences in $M[G]$), and
$\En$ is still a countable network for $\Cee$. 

On the other hand, the generic filter $G$ introduces an uncountable
strictly decreasing chain of open subsets of $Y(A)$. Thus,
$Y(A)^{M[G]}$ contains an uncountable free sequence, which contradicts
Theorem \ref{ciasnota}.
\end{pf}

\section{Some results concerning classes $\elsig(\loe n)$}

\begin{prop}\label{le1} Let $\ntr$ and assume $X$ is a space which has a
disjoint family of open sets $\setof{U_\al}{\al<\omega_1}$, and for
each $\alpha<\omega_1$ there is a closed set $Y_\alpha\subset
U_\alpha$ such that $Y_\al\notin\elsig(\loe n)$. Then
$X\notin\elsig(\loe n+1)$.
\end{prop}

\begin{pf} Suppose $X\in\elsig(\loe n+1)$ and fix a cover
$\Cee\subs\dpower X{\loe n+1}$ which has a countable network
$\En$. For each $\al<\omega_1$ the collection $\setof{Y_\alpha\cap
C}{C\in\Cee}$ is a cover of $Y_\al$ with a countable network in
$Y_\alpha$. Since $Y_\al\notin\elsig(\loe n)$, there exists 
$C_\alpha\in \Cee$ such that $C_\alpha\subs Y_\alpha$. Choose
$N_\alpha\in\En$ so that $C_\alpha\subset N_\alpha\subset
U_\alpha$. Then $N_\alpha\ne N_\beta$ for $\alpha\ne\beta$, a
contradiction. \end{pf}

\begin{prop}\label{le2} Let $\ntr$ and assume
$\setof{X_\xi}{\xi<\kappa}\subs \elsig(\loe n)$ is a family of compact
spaces and $\kappa\loe\cont$. Let $X$ be the one-point
compactification of\/ $\bigoplus_{\xi<\kappa}X_\xi$. Then $X\in
\elsig(\loe n+1)$. If $X_\xi\in \elsig(n)$ for uncountably many $\xi$,
then $X\in\elsig(n+1)$. 
\end{prop}

\begin{pf} For each $X_\xi$ fix a cover
$\Cee_\xi\subs\dpower{X_\xi}{\loe n}$ with a countable network
$\En_\xi$ in $X_\xi$. We assume that $X_\xi\cap X_\eta=\emptyset$
whenever $\xi\ne\eta$ and that
$X=\sn\infty\cup\bigcup_{\xi<\kappa}X_\xi$. 

Define $\Cee=\setof{C\cup\sn\infty}{C\in C_\xi,\;\xi<\kappa}$. We
will show that $\Cee$ has a countable network in $X$. 

Let $\En_\xi=\setof{N^\xi_k}{k<\omega}$ and fix a countable family
$\Aaa\subs\power\kappa$ which separates finite sets (here we use the
fact that $\kappa\loe\cont$). Let $A^*=\bigcup_{\xi\in
A}X_\xi$. Define
$$\En=\bigsetof{\sn\infty\cup\Bigl(A^*\cap\bigcup_{\xi<\kappa}N^\xi_k\Bigr)}
{k<\omega,\;A\in\Aaa}.$$

Then $\En$ is a countable family. We claim that $\En$ is a network for
$\Cee$. Fix $C=C_0\cup\sn\infty$, where $C_0\in\Cee_\xi$. Fix an open
set $U\subs X$ such that $C\subs U$. Then $U$ is a neighborhood of
$\infty$, so the set $F=\setof{\eta<\kappa}{X_\eta\not\subs U}$ is
finite. Find $A\in\Aaa$ such that $\xi\in A$ and
$(F\setminus\sn\xi)\cap A=\emptyset$. Find $k<\omega$ such that
$C_0\subs N^\xi_k\subs U\cap X_\xi$. Let $M=\sn\infty\cup
(A^*\cap\bigcup_{\eta<\kappa}N^\eta_k)$. Then $M\in \En$ and $C\subs
M\subs U$.

The second statement follows from Proposition \ref{le1}.
\end{pf}

\begin{wn} All classes $\elsig(n)$, for $n<\omega$, as well as
$\elsig(\omega)$, restricted to compact spaces, are nonempty. In
fact, for every $\ntr$ there exists a scattered compact space $X_n$ of
height $n$ and of cardinality $\omega_1$, such that
$X_n\in\elsig(n+1)$. \hfill$\Box$ \end{wn}

By Proposition \ref{ls-products},

\begin{prop}\label{p3} Assume $X\in \elsig(\loe n)$ and $Y\in
\elsig(\loe k)$, where $n,k<\omega$. Then $X\times Y\in\elsig(\loe
nk)$. \hfill$\Box$\end{prop}


Proposition \ref{le2} implies that $A(\kappa)\in\elsig(2)$ if
$\omega<\kappa\loe \cont$. Below we show that
$A(\omega_1)^n\in\elsig(n+1)$, but $A(\omega_1)\times
A(\omega_2)\notin\elsig(3)$ even if $\omega_2\loe\cont$.

\begin{tw} $A(\omega_1)^n\in\elsig(n+1)$. \end{tw}

\begin{pf} Let $A(\omega_1)=\omega_1+1$, where all points of $\omega_1$ are
isolated. The proof is by induction on $n$. The case $n=1$
is proved above. Suppose that $n>1$. Since
$A(\omega_1)^n$ contains $\omega_1$ disjoint clopen copies of
$A(\omega_1)^{n-1}$, $A(\omega_1)^n\not\in\elsig(<n+1)$ by Proposition~\ref{le1}.

For each permutation $p:n\rightarrow n$, consider the subset $X_p$ of
$A(\omega_1)^n$ defined by  
$$
X_p=\{(x_i)_{i<n}:x_{p(i)}\leq x_{p(i+1)}\text{ for all }i<n-1\}.
$$

Then $A(\omega_1)^n=\bigcup\{X_p:p\in{}^nn\text{ is a
permutation}\}$. By Proposition \ref{ls-invariance}, it suffices to prove that
each $X_p\in \Elsig(\leq n+1)$. Since all $X_p$ are homeomorphic, it
suffices to prove that $X_{id}\in \elsig(\leq n+1)$. We first define a
cover by $n+1$ element sets:

For each $\overline{\alpha}=(\alpha_0,\alpha_1,...,\alpha_{n-1})\in
X_{id}\cap \omega_1^n$ and for each $k<n$ let $\overline{\alpha}^k$ be
the element of $X_{id}$ obtained by changing the last $n-k$ coordinates
of $\overline{\alpha}$ to $\omega_1$. So,
$\overline{\alpha}^n=\overline{\alpha}$ and
$\overline{\alpha}^0=(\omega_1,\omega_1,...,\omega_1)$. Let 
$$
C_{\overline \alpha}=\{\overline{\alpha}^k:k\leq n,\ \overline{\alpha}\in X_p\cap \omega_1 ^n\}.
$$
These sets form our cover of $X_{id}$ by $n+1$-element sets. 

We now define a countable network for this family in $X_{id}$. The elements
of this network will be constructed from three countable families of
sets:

\noindent Family 1: We assume, by induction, that the collection of
similarly defined sets in $A(\omega_1)^{n-1}$ has a countable
network (note that for $n=1$ this follows from
\ref{exs-basic}.3). Since $A(\omega_1)^{n-1}$ is naturally identified with the
subspace  $Y=A(\omega_1)^{n-1}\times \{\omega_1\}$, there is a
countable network for  the sets $C_{\overline{\alpha}}\cap Y$
contained in $Y\cap X_{id}$. Call this countable family ${\mathcal
N}_Y$.  Notice that for each $\overline{\alpha}\in X_{id}$, the only
point of $C_{\overline{\alpha}}$ not contained in $Y$ is $\overline{\alpha}$. 

\noindent Family 2: Fix a countable family $\mathcal F\subseteq {}^{\omega_1}\omega_1$ of
functions such that
$\{(\alpha,\beta):\beta\leq\alpha\}=\bigcup{\mathcal F}$. For a
sequence $\overline{f}=(f_i:i<n-1)$ in ${\mathcal F}$, let
$N_{\overline{f}}\subseteq \omega_1^n$ be the set of
$(\alpha_i)_{i<n}$ such that $\alpha_i=f_i(\alpha_{i+1})$ for each
$i<n-1$. Notice that each $N_{\overline{f}}\subseteq X_{id}$. 

\noindent Family 3: Let ${\mathcal N}_0$ be a countable family of
subsets of $\omega_1$ that separates points from finite sets.

Our network will consist of sets of the following form:

$$
N \cup \bigr(N_{\overline{f}}\cap \prod\limits_{i<n} N_i\bigr)
$$
Where $N\in {\mathcal N}_Y$, $\overline f\in {\mathcal F}^{n-1}$ and
each $N_i\in {\mathcal N}_0$. 

To verify that this family is a network with respect to the cover by
the sets $C_{\overline{\alpha}}$, fix $\overline{\alpha}\in X_{\id}$
and an open set $U\supseteq C_{\overline{\alpha}}$. By construction,
there is an $N\in {\mathcal N}_Y$, such that $C_{\overline
\alpha}\setminus\{\overline{\alpha}\}\subseteq N\subseteq U$. Since
$\overline{\alpha}_0=(\omega_1,\omega_1,...,\omega_1)\in U$, we may
fix finite sets $F_i\subseteq \omega_1$ for $i<n$ so that
$$
V(0)=\prod\limits_{i<n} (\omega_1\setminus F_i)\subseteq U.
$$
Similarly, for each $0<k\leq n$, we may fix a basic open set $V(k)$
containing $\overline{\alpha}_k$. Without loss of generality, we may
assume that  
$$
V(k)=\{\alpha_0\}\times\dots\times\{\alpha_{k-1}\}\times\prod\limits_{k\leq i<n}
(\omega_1\setminus F_i)\subseteq U.
$$
Moreover, we may assume that $\alpha_i\in F_i$ for each $i<n$.

For each $i<n$ fix $N_i\in{\mathcal N}_0$ such that $N_i\cap
F_i=\{\alpha_i\}$ (that is, $N_i$ separates the point $\alpha_i$ from
the finite set $F_i\setminus \{\alpha_i\}$).

Next, for each $i<n-1$, fix $f_i\in {\mathcal F}$ such that
$f_i(\alpha_i+1)=\alpha_i$. Clearly,
$$
\overline{\alpha}\in N_{(f_i)}\cap \prod\limits_{i<n} N_i,
$$ 
so it suffices to prove that $N_{(f_i)}\cap \prod_{i<n} N_i\subseteq
U$. So suppose not, and fix
$\overline{\beta}=(\beta_0,\beta_1,...,\beta_{n-1})\in (N_{(f_i)}\cap
\prod_{i<n} N_i)\setminus U$. Fix $k$ maximal such that
$\beta_k=\alpha_k$ (if there is no such $k$, then $\beta_k\not\in F_k$
for every $k$ and hence $\overline{\beta}\in V(0)\subseteq U$). Then
$\beta_i=\alpha_i$ for all $i\leq k$, since $\overline\beta\in
N_{(f_i)}$, and $\beta_i\not\in F_i$ for all $i>k$. Thus,
$\overline{\beta}\in V(k)\subseteq U$, a contradiction.
\end{pf}

Let us now prove that $A(\omega_2)\times A(\omega_2)\notin \elsig(3)$.

\smallskip
Given two families $\Cal C_1$, $\Cal C_2$ of subsets of a set $X$, denote
$\Cal C_1\wedge \Cal C_2=\{\,C_1\cap C_2: C_1\in\Cal C_1,\ C_2\in\Cal C_2\,\}$.

\begin{lm} Let\/ $\Cal C_1$, $\Cal C_2$ be covers of a space $X$ by
finite sets, and $\Cal N_1$, $\Cal N_2$ networks with respect to 
$\Cal C_1$ and $\Cal C_2$. Then $\Cal N_1\wedge\Cal N_2$ is a network
with respect to $\Cal C_1\wedge\Cal C_2$.
\end{lm} 

\begin{pf}
Let $C_1=\{u_1,\dots,u_k, s_1,\dots, s_m\}\in \Cal C_1$,
$C_2=\{u_1,\dots, u_k, t_1,\dots,t_l\}\in\Cal C_2$, and let $U$ be a
neighborhood of $C=C_1\cap C_2=\{u_1,\dots, u_k\}$ in $X$.
Fix disjoint neighborhoods $U_{u_1}, \dots, U_{t_l}$ of the points in 
$C_1\cup C_2$ so that $U_{u_1}\cup\dots\cup U_{u_k}\subset U$.
Fix $N_1\in\Cal N_1$ and $N_2\in\Cal N_2$ so that 
$$
C_1\subset N_1\subset\bigcup\{\,U_{u_i}:1\loe i\loe k\,\}\cup\bigcup\{\,U_{s_i}:1\loe i\loe m\,\}
$$
and
$$
C_2\subset N_2\subset\bigcup\{\,U_{u_i}:1\loe i\loe k\,\}\cup\bigcup\{\,U_{t_i}:1\loe i\loe l\,\}
$$
Then $N=N_1\cap N_2$ is an element of $\Cal N_1\wedge\Cal N_2$ such
that $C_1\cap C_2\subset N\subset \bigcup\{\,U_{u_i}:1\loe i\loe k\,\}\subset U$.
\end{pf}

\smallskip
\begin{tw} $A(\omega_2)^2\notin \elsig(3)$. \end{tw}

\begin{pf} If $\omega_2> 2^\omega$, then $A(\omega_2)$ is not in
$\elsig(3)$, because every space in $\elsig(n)$ has cardinality 
$\loe2^\omega$.

So assume $\omega_2\loe 2^\omega$. Let $A(\omega_2)=\omega_2+1$, where
all points of $\omega_2$ are isolated. Let $C_0$ be the cover
of $A(\omega_2)^2$ by the sets of the form 
$C_{\alpha\beta}=\{(\alpha,\beta), (\alpha,\omega_2), (\omega_2,\beta),
(\omega_2,\omega_2)\}$, $\alpha,\beta <\omega_2$. Then there is a
countable network $\Cal N_0$ with respect to $\Cal C_0$ (see Example
\ref{exs-basic}(3) and equivalence of (a) and (b) in Proposition
\ref{p2-0}).

Now suppose that there exists a cover $\Cal C$ of $A(\omega_2)^2$ with
at most 3-point sets and a countable network $\Cal N$ with respect to
$\Cal C$. Replacing $\Cal C$ with $\Cal C\wedge \Cal C_0$ and $\Cal N$
with $\Cal N\wedge \Cal N_0$ if necessary, we may assume that every
element of $\Cal C$ is contained in a set of the form $C_{\alpha\beta}$.

We will say that an element $C$ of $\Cal C$ is of type 1 if for some
$\alpha, \beta$, $C\subset \{(\alpha,\beta), (\omega_2,\beta),
(\omega_2,\omega_2)\}$, is of type 2 if $C\subset \{(\alpha,\beta),
(\alpha,\omega_2),(\omega_2,\omega_2)\}$, and of type 3 if
$C\subset \{(\alpha,\beta), (\alpha,\omega_2), (\omega_2,\beta)\}$.
Note that the elements of $\Cal C$  of the three types cover
$\omega_2\times \omega_2$. Furthermore, there are at most countably
many elements of type 3. Indeed, otherwise the union $P$ of all
elements of $\Cal C$ of type 3 would be a Lindel\"of $\Sigma$-space, which is
impossible, because one of the three sets 
$$
P\cap (A(\omega_2)\times\{\omega_2\}),
\quad P\cap\bigl(\{\omega_2\}\times A(\omega_2)\bigr), \quad
\{\,(\alpha,\beta)\in P: (\alpha,\omega_2)\notin
P,\ (\omega_2,\beta)\notin P\,\}
$$
would be an uncountable closed discrete subspace of $P$. Removing from
$A(\omega_2)$ all elements of $\omega_2$ that occur as coordinates of
points of $P$, and taking intersections of $\Cal C$ the square of the
remaining set, we will obtain a cover without elements of type 3
and the intersection of $\Cal N$ with this square is a countable
network with respect to this new cover. Thus, we may assume that all
elements of $\Cal C$ are of type 1 or 2.

Let $(\alpha,\beta)\in\omega_2\times\omega_2$; choose an element C of
$\Cal C$ such that $(\alpha,\beta)\in C$. If $C$ is of type 1, then
$V=(A(\omega_2)\times\{\beta\})\cup (A(\omega_2) \setminus\{\alpha\})^2$
is a neighborhood of $C$, so there is an $N\in \Cal N$ such
that $C\subset N\subset V$. Obviously, 
$N\cap (\{\alpha\}\times \omega_2)=\{(\alpha,\beta)\}$.
Similarly, if $C$ is of type 2, there is an $N\in\Cal N$ such that
$N\cap (\omega_2\times\{\beta\})=\{(\alpha,\beta)\}$.
Thus, the sets
$$
B_1=\{\,(\alpha,\beta)\in\omega_2\times\omega_2: N\cap (\{\alpha\}\times
\omega_2)=\{(\alpha,\beta)\}\text{ for some $N\in\Cal N$}\,\}
$$
and 
$$
B_2=\{\,(\alpha,\beta)\in\omega_2\times\omega_2:N\cap (\omega_2\times\{\beta\})=\{(\alpha,\beta)\}
\text{ for some $N\in\Cal N$}\,\}
$$
cover $\omega_2\times\omega_2$. Note that for any given $\alpha\in
\omega_2$ and $N\in\Cal N$ there is at most one $\beta\in\omega_2$
such that $N\cap (\{\alpha\}\times \omega_2)=\{(\alpha,\beta)\}$, so
for every $\alpha\in\omega_2$, the set 
$B_1\cap (\{\alpha\}\times\omega_2)$ is at most countable. Similarly,
for every $\beta\in\omega_2$ the set
$B_2\cap(\omega_2\times\{\beta\})$ is at most countable. The existence
of such a pair of sets covering $\omega_2\times\omega_2$ contradicts a
theorem of Kuratowski \cite{kura}.
\end{pf}

Note that the same proof shows that $A(\omega_1)\times A(\omega_2) \not\in L\Sigma(3)$. 

\vskip 6pt

It appears natural to expect that if $X\in \elsig(n)$ for some $n>1$,
then the sequence $n_k$ such that $X^k\in \elsig(n_k)$ should increase
to infinity; hence the following question:

\begin{question} \label{countable-power} Suppose $X^\omega\in
\elsig(<\omega)$. Must $X$ have a countable network?\end{question}

The results in this section allow to obtain a consistently positive
answer to this question.

\medskip
\begin{prop}\label{finite-decompo} Suppose $X\in \elsig(<\omega)$. Then
there are subspaces $X_n\subset X$, $n\in\omega$, such that
$X_n\in\elsig(\loe n)$, and $X=\bigcup\{\,X_n:n\in\omega\,\}$.
\end{prop}

\begin{pf} Let $p\:M\to X$ be a finite-valued usc mapping from a
second-countable space $M$ such that $p(M)=X$. For every $n\in\omega$,
put $M_n=\{\,m\in M: |p(m)|\le n\,\}$ and $X_n=p(M_n)$.
\end{pf}

\begin{prop}\label{powers-1} Assume $X$ is a space such that
$X^{\omega}\in\elsigfin$. Then for some $n\in\omega$, 
$X^\omega\in\elsig(n)$.
\end{prop}

\begin{pf} For every $k\in\omega$, let 
$\pi_k\:(X^\omega)^\omega\to X^\omega$ be the projection to the $k$th
factor. Since $(X^\omega)^\omega$ is homeomorphic to $X^\omega$, 
we have $(X^\omega)^\omega=\bigcup\{\,X_k:k\in\omega\,\}$ where
$X_k\in\elsig(\loe k)$. Obviously, $\pi_n(X_n)=X^\omega$ for some
$n\in\omega$. Then $X^\omega\in \elsig(\loe n)$.
\end{pf}

\begin{tw} If $X^\omega\in\elsig(<\omega)$, then $X^\omega$ is hereditarily
separable.\end{tw}  

\begin{pf}  Let us first show that $X^\omega$ has no uncountable
discrete subspaces. By Lemma~\ref{powers-1}, $X^\omega\in \elsig(n)$
for some $n\in \omega$. Let $D$ be a discrete subspace of $X$, and
let $F$ be its closure in $X$. Since $F^\omega$ is a closed subspace
of $X^\omega$,  $F^\omega\in \elsig(\loe n)$. If $D$ were
uncountable, this would be impossible, because from Proposition~\ref{le1} by
a simple inductive argument would follow $F^n\notin\elsig(\loe n)$.
Thus, every discrete subspace of $X$ is countable.

Since $X^\omega$ is homeomorphic to its square, it follows that
$X^\omega$ is hereditarily separable or hereditarily Lindel\"of
\cite[Theorem 1]{Zenor}. If $X^\omega$ is hereditarily Lindel\"of,
then it has $G_\delta$-diagonal, and since it is a Lindel\"of
$\Sigma$-space, it must have countable network \cite{ava-mono}.
\end{pf}

Thus, if the answer to Question~\ref{countable-power} is negative,
then there is a strong $S$-space. Since $\MAone$ implies that
there are no strong $S$-spaces \cite{Kunen}, we get

\begin{wn} If $\MAone$ holds, and $X^\omega\in\elsig(<\omega)$,
then $X$ has countable network.\end{wn}

\section{Dense metrizable subspaces}

Tka\v cenko asks in \cite{tkacenko} whether every compact
space in $\elsigO$ has a dense metrizable subspace. One cannot hope
for a completely metrizable dense subspace because, for example, all
metrizable subspaces of the double arrow space (which is in
$\kelsig(\loe2)$) are countable.

We give a partial answer to Tka\v cenko's question, namely, we show
that every space in $\kelsig(<\omega)$ has a dense metrizable
subspace. On the other hand, under $\neg \MAone$, there exists a
non-separable ccc space $X\in\kelsigO$ which is even metrizably
fibered; see Theorem 3.5 in the survey article of Todor\v cevi\'c
\cite{Todorcevic2000}. Observe that $X$ cannot have a dense metrizable
subspace, because such a subspace would be ccc and therefore
separable. This gives a consistent negative answer to Tka\v cenko's
question (see Problems 3.1 and 3.5 in \cite{tkacenko}).

Todor\v cevi\'c proved in \cite{Todorcevic1999} that every Rosenthal
compact space has a dense metrizable subspace. The key fact in the
proof is the absoluteness of the class of Rosenthal compacta with
respect to forcing extensions. We use the same idea, proving that the
class $\kelsig(<\omega)$ is absolute.

The following result is due to H.E. White, Jr. \cite{white}.

\begin{prop}[White Jr. \cite{white}]\label{klucz1} Every first
countable Hausdorff space with a $\sig$-disjoint \ppi-base contains a
dense metrizable subspace.
\end{prop}

Let $X$ be a compact space and let $\lat=\closed(X)$, the collection
of all closed subsets of $X$ or let $\lat$ be a sublattice of
$\closed(X)$ which is at the same time a closed base (we shall say
that $\lat$ is a {\em basic lattice} for $X$). Observe that every
point of $X$ corresponds to an ultrafilter (~=~a maximal filter) in
$\lat$. More precisely, the compact space $X$ can be recovered from
$\lat$ as the space $\ult(\lat)$ of all ultrafilters in $\lat$; the
topology is generated by closed sets of the form
$a^+=\setof{p\in\ult(\lat)}{a\in p}$, where $a\in \lat$. This gives an
idea of how to interpret a compact space in a forcing extension. Namely,
if $X$ is a compact space in a ground model (a countable transitive
model of a ``large enough" part of ZFC) $M$ and if $G$ is a
$\poset$-generic filter over $M$, where $\poset\in M$ is a fixed
forcing notion, then we can define $X[G]$, the {\em interpretation of}
$X$ in $M[G]$, as $\ult(\lat)$ computed in $M[G]$, where
$\lat=\closed(X)$ (in fact, $\lat$ can be any basic lattice for
$X$). In $M[G]$, $\lat$ is the same algebraic object, but
$\ult^{M[G]}(\lat)$ is usually different from $\ult^M(\lat)$. It can
be proved that the definition of $X[G]$ does not depend on the choice
of $\lat$ (as long as $\lat$ is a basic lattice for $X$).

Let us remark that compact spaces in
forcing extensions and absoluteness were already studied by Bandlow
\cite{Bandlow}.

A lattice $\seq{\lat,\join,\meet,0,1}$ is {\em normal} if it is
distributive and for every $a,b\in \lat$ such that $a\meet b=0$ there
exist $a',b'\in\lat$ such that 
$$a\meet b'=a'\meet b=0\quad\text{ and }\quad a'\join b'=1.$$

Every basic lattice for a compact space is normal.

The following result is due to Todor\v cevi\'c \cite{Todorcevic1999},
although it is not stated explicitly in \cite{Todorcevic1999}. For
the proof see Lemma 4 of \cite{Todorcevic1999}.

\begin{prop} \label{klucz2} Assume $X$ is a compact space such that for every $\ro(X)$-generic
filter $G$ over some ground model containing $X$, the extension $X[G]$
has countable tightness. Then $X$ has a $\sig$-disjoint \ppi-base.
\end{prop}

The following result appears in Tkachuk \cite{tkachuk}. For the sake
of completeness we give a proof, which is different from the one in
\cite{tkachuk}.

\begin{lm}\label{klucz3} Assume $K\in\elsigO$ is compact. Then $K$ has
a dense set of $G_\delta$-points. \end{lm}



\begin{pf} It suffices to show that $K$ contains at least
one $G_\delta$-point, since the class $\elsigO$ is stable under closed
subsets. Let $\Cee$ be a compact cover of $K$ consisting of metrizable
sets and let $\En=\setof{N_n}{\ntr}$ be a countable network for $\Cee$
which consists of closed sets. Choose a sequence of open sets
$\setof{U_n}{\ntr}$ such that $\cl U_{n+1}\subs U_n$ and either
$U_n\subs N_n$ or $U_n\cap N_n=\emptyset$. Let 
$F=\bigcap_{\ntr}U_n$. Then $F$ is a nonempty closed $G_\delta$-set. Choose $C\in \Cee$ such that $F\cap C\nnempty$. If
$F\not\subs C$, then there is an $N\in\En$ such that $C\subs N$ and 
$F\not\subs N$, which is impossible by the construction of $F$. Thus,
$F\subs C$. It follows that $F$ is a closed $G_\delta$-set contained in a
metrizable subspace $C$ of $K$. Hence, every point of $F$ is $G_\delta$
in $F$, and therefore also $G_\delta$ in $K$.
\end{pf}

The rest of this section is devoted to the proof that the class
$\kelsig(<\omega)$ is absolute.  

Assume $\map\phi XY$ is a usc compact-valued function with nonempty
images of points, and let $\K =\closed(X)$, $\lat=\closed(Y)$. Define
$\map h\lat\K$ by setting

$$h(a)=\setof{x\in X}{\phi(x)\cap a\nnempty}.$$

Note that $h$ is well defined by the upper semicontinuity of
$\phi$. We will call $h$ the {\em lattice map associated} to 
$\phi$. Assume further that $X$ is compact and $Y=\bigcup_{x\in
X}\phi(x)$. Observe that $h$ has the following properties:

\begin{enumerate}
	\item[(1)] $h(a_1\cup a_2)=h(a_1)\cup h(a_2)$.
	\item[(2)] $h(\emptyset)=\emptyset$ and $h(Y)=X$.
	\item[(3)] $h(a)\nnempty$ whenever $a\nnempty$.
	\item[(4)] If $h(a)\cap b=\emptyset$ then there exists $c$
such that $a\cap c=\emptyset$, $b\subs h(c)$ and $b\cap
h(a')=\emptyset$ whenever $c\cap a'=\emptyset$.
\end{enumerate}

Properties (1)--(3) are obvious. To see (4), take $c=\bigcup_{x\in
b}\phi(x)$ and note that by usc, $c$ is a compact set, and therefore
is in $\lat$.

It turns out that it is possible to reconstruct $\phi$ from a
map $\map h\lat\K$ satisfying (1)--(4). More precisely:

\begin{lm}\label{normalne} Assume $\lat$, $\K$ are normal lattices and
$\map h\lat\K$ is a map satisfying conditions (1)--(4) above. Then
there exists a compact-valued usc map
$\map\phi{\ult(\K)}{\ult(\lat)}$ such that
\begin{equation}h(a)^+=\setof{p\in\ult(\K)}{\phi(p)\cap
a^+\nnempty}.\tag{$*$}\end{equation}
\end{lm} 

\begin{pf} Fix $p\in\ult(\K)$ and define

$$I(p)=\setof{a\in\lat}{h(a)\notin p}.$$

Observe that $I(p)=\inv hp$ is an ideal by (1), (2). Denote by ``$a\ll
b$" the relation ``$a^+\subs\Int b^+$", which by normality can be
defined algebraically as ``there exists $c$ such that $a\meet c=0$ and
$c\join b=1$". We claim that $I(p)$ is $\ll$-directed, i.e. for every
$a\in I(p)$ there is $a'\in I(p)$ with $a\ll a'$. 

Fix $a\in I(p)$. Find $b\in p$ with $h(a)\meet b=0$. Using (4) find
$c$ such that $a\meet c=0$, $b\loe h(c)$ and $b\meet h(a')=0$ whenever
$c\meet a'=0$. By normality, there are $a',c'$ such that $a'\join
c'=1$ and $a\meet c'=0=a'\meet c$. It follows that $a\ll a'$ and
$h(a')\notin p$ because $h(a')\meet b=0$. Thus $a'\in I(p)$.

The property that $I(p)$ is $\ll$-directed ensures that the set
$$
u(p)=\bigcup\setof{a^+}{h(a)\notin p}
$$ is open. Define
$$
\phi(p)=\ult(\lat)\setminus u(p).
$$

We need to check that ($*$) holds and that $\phi$ is usc.

Fix $p\in h(a)^+$ and suppose $\phi(p)\cap a^+=\emptyset$. Then
$a^+\subs u(p)$ which, by the fact that $I(p)$ is a $\ll$-directed
ideal, implies that $a\in I(p)$, a contradiction. Conversely, if
$\phi(p)\cap a^+\nnempty$, then certainly $a\notin I(p)$, which means
that $p\in h(a)^+$. This shows ($*$).

To see that $\phi$ is usc, fix an open set $u\subs \ult(\lat)$ and fix
$p$ such that $\phi(p)\subs u$. Using the fact that
$\setof{a^+}{a\in\lat}$ is a basic lattice for $\ult(\lat)$, find
$c\in \lat$ such that $u\cup c^+=\ult(\lat)$ and $\phi(p)\cap
c^+=\emptyset$. Then 
$$
p\in \setof{x\in\ult(\K)}{\phi(x)\cap
c^+=\emptyset}\subs\setof{x\in\ult(\K)}{\phi(x)\subs u},
$$
and the set in the middle is open by ($*$). It follows that $\phi$ is usc.

This completes the proof. \end{pf}

\begin{lm}\label{WuEf} Assume $\lat$, $\K$ are basic lattices for
compact spaces $Y$ and $X$ respectively and $\map h\lat\K$ is a map
associated to a compact-valued function $\map\phi XY$. Define

$$T_h=\setof{s\in\fin\lat}{(\forall\;a_1,a_2\in s)\;a_1\cap
a_2=\emptyset \text{ and }\bigcap_{a\in s}h(a)\nnempty}.$$

Endow $T_h$ with a strict partial order $<$ defined by $s<t$ iff
\begin{enumerate}
	\item[(i)] every $a\in s$ is below some $a'\in t$,
	\item[(ii)] for every $b\in t$ there is $a\in s$ with $a\subs b$ and
	\item[(iii)] for some $b\in t$ there are two distinct (and
therefore disjoint) sets $a_1,a_2\in s$ such that $a_1\cup a_2\subs
b$.
\end{enumerate}
Then $\phi$ is finite-valued if and only if $\pair{T_h}{<}$ is well-founded.

\end{lm}

\begin{pf} Assume that $s_1>s_2>\dots$ is an infinite decreasing
sequence in $\pair{T_h}<$. 
Observe that $\bigcap_{a\in s_n}h(a)\sups \bigcap_{a\in s_{n+1}}h(a)$
for every $\ntr$ and therefore by the compactness of $X$, there is
$x\in X$ such that $x\in  \bigcap_{a\in s_n}h(a)$ for every
$\ntr$. This means that $\phi(x)\cap a\nnempty$ whenever $a\in s_n$
and $\ntr$.

Let $F=\phi(x)\cap\bigcap_{\ntr}\bigcup s_n$. By compactness and by
the definition of $<$, $F$ is an infinite closed subset of $Y$, which
shows that $\phi(x)$ is infinite.

Assume now that $\phi(x)$ is an infinite set for some $x\in
X$. Construct inductively a sequence $s_1>s_2>\dots$ in $T_h$ such
that $a\cap \phi(x)\nnempty$ for every $a\in s_n$ and
$\phi(x)\subs\bigcup s_n$. Since $\phi(x)$ is infinite, $a\cap
\phi(x)$ is infinite for some $a\in s_n$, and therefore it is always
possible to find $s_{n+1}<s_n$ satisfying the above condition. This
shows that $\pair{T_h}<$ is not well-founded.
\end{pf}

\begin{tw} The class $\kelsig(<\omega)$ is absolute. That is, if $M$ is
a transitive model of (a large enough part of) ZFC such that $X\in M$
and $M\models X\in\kelsig(<\omega)$ then, setting\/ $\lat=\closed(X)$
(defined in $M$), for every ZFC model $N$ which extends $M$, we have
that $$N\models \;\ult(\lat)\in\kelsig(<\omega).$$
\end{tw}

\begin{pf} We work in $M$:
Let $T$ be a compact metric space and let $\map\phi TX$ be a
finite-valued usc function such that $X=\bigcup_{t\in T}\phi(t)$. Let
$\lat=\closed(X)$ and $\K=\closed(T)$. Let $\map h\lat\K$ be the
associated map. Then the poset $\pair{T_h}{<}$ is well founded by
Lemma \ref{WuEf}.

In $N$, $h$ is a map of normal lattices satisfying conditions (1)--(4)
(these conditions are clearly absolute), and $\pair {T_h}<$ is the same
object as in $M$. The property of being well-founded is
absolute. Thus, by Lemmas \ref{normalne} and \ref{WuEf}, we deduce
that $\ult^N(\lat)$ is an image of a finite-valued usc function
defined on $\ult^N(\K)$. It remains to observe that $\ult(\K)$ is
metrizable. The metrizability of $\ult(\K)$ is equivalent to the fact
that $\K$ is separated as a lattice by a countable family, that
is, there is a countable $N\subs \K$ such that for every $a,b\in \K$
with $a\meet b=0$ there are $a',b'\in N$ such that $a\loe a'$, $b\loe
b'$ and $a'\meet b'=0$. The last property is true in $M$ and remains
true in any ZFC model containing $M$.

\end{pf}

\begin{wn} Every space in $\kelsig(<\omega)$ contains a dense
metrizable subspace. \end{wn}

\begin{pf} By the above theorem combined with Lemma \ref{klucz3},
Propositions \ref{klucz1}, \ref{klucz2} and Corollary~\ref{tightness}. \end{pf}

Let us mention that, at least consistently, the class $\kelsigO$ is not
absolute. In fact, the example of Todor\v cevi\'c from
\cite{Todorcevic2000} (mentioned in the beginning of this section)
constructed under $\neg\MAone$, which is a metrizably fibered compact
space, must have uncountable tightness after forcing with its
regular-open algebra --- otherwise it would contain a dense metrizable
subspace.


\section{Results under $\MAone$}

We prove that $\MAone$ implies that each compact spaces of scattered
height $3$ and of size $\omega_1$ is in $\elsig(\loe 3)$. We first
state and prove a general lemma about strongly almost disjoint
families on $\omega_1$. By a strongly almost disjoint family of sets
we mean any collection of infinite sets with pairwise intersections
finite. The family may contain both countable and uncountable sets. 

\begin{lm}\label{main_lemma} Suppose that
$\setof{A_\alpha}{\alpha<\omega_1}$ is a strongly almost disjoint
family of subsets of\/ $\omega_1$. Suppose also that
$\setof{p_\alpha}{\alpha<\omega_1}$ is a sequence of pairwise disjoint
finite subsets of $\omega_1$. Then there are $\alpha<\beta$ such that
$p_\alpha\cap A_\beta=\emptyset$ and $p_\beta\cap
A_\alpha=\emptyset$. 
\end{lm}

\begin{pf} By passing to a subsequence, we may assume that there is
$n\in \omega$ such that $|p_\alpha|=n$ for all $\alpha$.  Let $M$ be a
countable elementary submodel containing everything relevant and let
$\gamma=M\cap \omega_1$.

\begin{claim}\label{sub} There are
$\setof{\alpha_i}{i<n+1}\subs\gamma$
and a $\beta>\gamma$ such that 
$$
(\bigcup_{i<n+1} p_{\alpha_i})\cap A_\beta = \emptyset
$$
\end{claim}

\begin{pf} If not, then for each $\beta>\gamma$ there are at most $n$
$\alpha$'s below $\gamma$ such that $p_\alpha\cap
A_\beta=\emptyset$. Thus there is a $\alpha_\beta<\gamma$ such that
$p_\alpha\cap A_\beta\not=\emptyset$ for each $\alpha\in M\setminus
\alpha_\beta$. Choose $\setof{\beta_i}{i<n+1}\subs M\setminus
\gamma$. Choose an index $\alpha\in M$ above
$\setof{\alpha_{\beta_i}}{i<n+1}$ such that all ordinals in $p_\alpha$
lie above the following finite set:
$$
 \bigcup_{0\loe i<j<n+1} (A_{\beta_i}\cap A_{\beta_j})\cap M
$$
Then for $i<n+1$ we have that $p_\alpha\cap A_{\beta_i}\not=\emptyset$
and the sets $p_\alpha\cap A_{\beta_i}$ are disjoint. This contradicts
that $|p_\alpha|=n$. 
\end{pf}

To complete the proof of the main lemma, note that since
$A_{\alpha_i}\cap A_{\alpha_j}\subs M$ for each $i\not=j$, by the
pigeon-hole argument just presented in the proof of Claim \ref{sub},
we may conclude that $p_\beta\cap A_{\alpha_i}=\emptyset$ for some
$i<n+1$ (here we use that $p_\beta\cap M=\emptyset$ for
$\beta>\gamma$, and this can be easily arranged by going to a
subsequence).
\end{pf}

\begin{tw}\label{MA} $\MAone$ implies that all compact spaces of
cardinality $\omega_1$ and scattered of height $3$ are in $\elsig(\loe
3)$. 
\end{tw}

\begin{pf} Let $X$ be a compact space of scattered height $3$ and of
cardinality $\omega_1$. Without loss of generality, $X$ is of the form
$I\cup D\cup \{\infty\}$ where $I=\omega_1$ is the set of isolated
points and $D=\setof{d_\alpha}{\alpha<\omega_1}$ the set of isolated
points in $X\setminus \omega_1$ (the cases where either $I$  or $D$ is
countable do not require any extra set-theoretic assumptions and the
space is in $\elsig(\loe 2)$).

For each $\alpha<\omega_1$ let $U_\alpha$ be a clopen neighborhood of
the point $d_\alpha$ such that $U_\alpha\setminus I=\{d_\alpha\}$. Let
$a_\alpha=U_\alpha\cap \omega_1$. Then
$\setof{a_\alpha}{\alpha<\omega_1}$ is a strongly almost disjoint family of
subsets of $\omega_1$, i.e., $a_\alpha\cap a_\beta$ is finite for each
$\alpha\ne\beta$. Some, but not necessarily all, of the sets
$a_\alpha$ may be uncountable. We may assume that 

\begin{enumerate}
\item[(a)] For each finite $F\subs \omega_1$,
$\omega_1\setminus \bigcup_{\alpha\in F} a_\alpha$ is uncountable.
\end{enumerate}
If not, it easy to see that $X$ is the sum of a countable subspace and a
compact subspace with finitely many non-isolated points, and therefore
(without any extra set-theoretic assumptions) $X$ is in $\elsig(\loe 2)$.
If we choose the neighborhoods $U_\alpha$ with some care, we may assure that 
\begin{enumerate}
\item[(b)] for each finite $x\subs \omega_1$ and each 
$\alpha<\omega_1$, $x$ is covered by $\bigcup_{\beta>\alpha} a_\beta$. 
\end{enumerate}
Therefore, $\infty$ has a local base consisting of sets of the form 
$$
\{\infty\} \cup (D\setminus F)\cup (\omega_1\setminus \bigcup_{\alpha\in F}a_\alpha)
$$
We may also make sure that 
\begin{enumerate}
\item[(c)] for each finite $x\subs \omega_1$ there are uncountably
many $\alpha$ such that $a_\alpha \cap x=\emptyset$.
\end{enumerate}
We now define a poset ${\mathbb P}$ consisting of all pairs $(p,F)$
where both $p,F\in [\omega_1]^{<\omega}$ with the property 
that $p\cap \bigcup_{\alpha\in F} a_\alpha =\emptyset$. We define
$(p,F)<(q,G)$ if $q\subs p$ and $G\subs F$.

\begin{claim} ${\mathbb P}$ has the ccc.
\end{claim}

\begin{pf} Let $\setof{(p_\alpha,F_\alpha)}{\alpha<\omega_1}$ be an
uncountable subset of ${\mathbb P}$. By a standard $\Delta$-system
argument, we may assume that the sets $p_\alpha$'s are pairwise disjoint and
likewise for the $F_\alpha$'s. Let $A_\alpha=\bigcup_{\beta\in
F_\alpha} a_\beta$.  By Lemma \ref{main_lemma}, we may find
$\alpha<\beta$ such that  $p_\alpha\cap A_\beta=\emptyset=p_\beta\cap
A_\alpha$. It follows that $(p_{\alpha}\cup p_\beta, F_{\alpha}\cup
F_\beta)\in {\mathbb P}$ is a common extension.
\end{pf}

Let $\mathbb P_\omega$ denote the finite support product of countably
many copies of ${\mathbb P}$. By $\MAone$, ${\mathbb P}_\omega$ has the ccc,
and we may fix a filter $G$ in ${\mathbb P}_\omega$ generic for the
following dense sets:

$$
E_{\alpha, G}=\{\,r\in {\mathbb P}_\omega: r(n)=(p,F) \text{ for some
$n\in\omega$, $p\ni\alpha$ and $F\supset G$}\,\}
$$ 
for $\alpha\in \omega_1$ and $G\subs \omega_1$ finite such that $\alpha\not\in \bigcup_{\beta\in G}a_\beta$.

One easily defines from the generic two countable families of subsets of $\omega_1$, 
$\setof{I_n}{n\in \omega}$ and $\setof{D_n}{n\in \omega}$ such that
\begin{enumerate}
\item[(d)] $\omega_1=\bigcup_n I_n=\bigcup_n D_n$

\item[(e)] For each $n$, $I_n\cap a_\alpha=\emptyset$ for each $\alpha\in D_n$

\item[(f)] For each $\alpha\in \omega_1$ and each finite $F\subs \omega_1$ such that $\alpha\not\in \bigcup_{\beta\in F}a_\beta$, there is an $n$ such that $\alpha\in I_n$ and $F\subs D_n$. 
\end{enumerate}

We need to define a similar family of sets before we give the cover
and network witnessing that  $X\in \elsig(\loe 3)$. To do this we let
$a_\alpha^\prime=a_\alpha\cap \bigcup_{\beta<\alpha}a_\beta$ for each
$\alpha$. Thus,

\begin{enumerate}
\item[(g)] $\setof{a_\alpha\setminus a_\alpha^\prime}{\alpha\in \omega_1}$ is a disjoint family 
and $\bigcup\setof{a_\alpha\setminus a_\alpha^\prime}{\alpha\in \omega_1}=\omega_1$. 
\end{enumerate}
Now we define another poset ${\mathbb Q}$ consisting of all pairs
$(p,F)$ where $p,F\in [\omega_1]^{<\omega}$ with the property
that $|p\cap a_\alpha^\prime|=1$ and $p\cap a_\alpha=p\cap
a_\alpha^\prime$ for each $\alpha\in F$. I.e., for each $\alpha\in F$,
$p$ intersects $a_\alpha$ at exactly one point, and that point is in
$a_\alpha^\prime$.

We take the same ordering given by $\sups$ on both coordinates.

\begin{claim} ${\mathbb Q}$ has the ccc.
\end{claim}

\begin{pf} Given an uncountable subset
$\setof{(q_\alpha,G_\alpha)}{\alpha\in \omega_1}$, we may assume that
the $q_\alpha$'s form a $\Delta$-system with root $r$. Let
$p_\alpha=q_\alpha\setminus r$. We may also assume that the
$G_\alpha$'s form a $\Delta$-system with root $R$. Let
$F_\alpha=G_\alpha\setminus R$. By going to a subsequence we may
assume that $p_\alpha\cap \bigcup\setof{a_\xi^\prime}{\xi\in
R}=\emptyset$ for each $\alpha$. Thus, for each $\xi\in R$,
$q_\alpha\cap a_\xi^\prime=r\cap a_\xi^\prime$. Thus, if we let
$A_\alpha=\bigcup_{\xi\in F_\alpha}a_\xi$ it follows that
$(q_\alpha,G_\alpha)$ is compatible with some $(q_\beta,G_\beta)$ if
and only if $p_\alpha\cap A_\beta=\emptyset=p_\beta\cap A_\alpha$. The
existence of such a pair is given by Lemma \ref{main_lemma}.
\end{pf}

Again we take the finite support product of countably many copies of
${\mathbb Q}$ and denote it by ${\mathbb Q}_\omega$. Taking dense sets
defined similarly as above, $\MAone$ gives us two countable families
$\setof{J_n}{n\in \omega}$ and $\setof{E_n}{n\in \omega}$ such that
\begin{enumerate}
\item[(h)] $|J_n\cap a_\alpha|=|J_n\cap a_\alpha^\prime|=1$ for each
$n\in \omega$ and for each $\alpha\in E_n$.
\item[(i)] For each $\beta\in \omega_1$ and for each finite $F\subs
\omega_1$ such that $\beta\in a_\alpha^\prime$ for each $\alpha\in F$,
there is an $n$ such that $\beta\in J_n$ and $F\subs E_n$.
\end{enumerate}

We are now ready to describe a cover $C$ of the space $X$ consisting
of 3-element sets and a countable family of subsets that is a network
for the cover.

For each $\alpha\in\omega_1$ and each $\beta\in a_\alpha\setminus
a_\alpha^\prime$ put $\{\beta,d_\alpha,\infty\}\in C$. Also put
$\{d_\alpha,\infty\}\in C$ in the case that
$a_\alpha=a_\alpha^\prime$. Clearly, by clause (g), this collection of
sets covers $X$. Now define a countable family of sets as follows: Let
$N_D$ be a countable family of subsets of $D$ that separates points
from finite sets. Then let
${\mathcal N}$ consist of all sets of the following forms:
\begin{enumerate}
\item $\{\infty\}\cup N\cup(I_n\cap J_m)$ for $N\in N_D$ and $n,m\in \omega$. 
\item $\{\infty\}\cup N\cup I_n$ for $N\in N_D$ and $n\in \omega$,
\item $\{\infty\}\cup N\cup J_n$ for $N\in N_D$ and $n\in \omega$, and
\item $\{\infty\}\cup N$, for $N\in N_D$.
\end{enumerate} 

To see that this countable set is a network at $C$, fix
$c=\{\beta,d_\alpha,\infty\}\in C$ and fix an open set $U\sups c$ (if
$c$ is of the form $\{d_\alpha,\infty\}$ the proof is easier). By
clause (b), there is a finite $H\subs \omega_1$ with  $\alpha\not\in H$
such that 
$$
V_H= \{\infty\}\cup(D\setminus H)\cup \bigl(\omega_1\setminus \bigcup_{\xi\in H}a_\xi\bigr)\subs U.
$$
Thus, $c\subs V_H\cup\{\beta\}\subs U$. 
Let $F=\setof{\xi\in H}{\beta\not\in a_\xi}$ and let $G=\setof{\xi\in
H}{\beta\in a_\xi}$. Assume that $F$ and $G$ are both not empty (in
the case that one is empty, the proof is similar). Note that since
$\beta\in a_\alpha\setminus a_\alpha^\prime$, it follows that
$\beta\in a_\xi^\prime$ for each $\xi \in G$. Fix $n$ such that
$\beta\in I_n$ and $F\subs D_n$, and fix $m$ such that $\beta\in J_m$
and $G\subs E_m$. Then it is straightforward to verify that $\beta\in
I_n\cap J_m\subs V_H\cup\{\beta\}$. Also fix $N\in N_D$ such that
$\alpha\in N$ and $N\cap H=\emptyset$. So, $c\subs \{\infty\}\cup
N\cup (I_n\cap J_m)\subs U$ as required. This completes the proof of
the theorem.
\end{pf}




We now show that some assumption is needed for Theorem \ref{MA}. Consider the following
$\clubsuit$-like principle: 

\begin{enumerate}
\item[$(*)$] There is an almost disjoint family
$\setof{a_\alpha}{\alpha\in \omega_1}$ of countable subsets of
$\omega_1$ such that for each countable family $X\subs
[\omega_1]^{\omega_1}$ there are uncountably many $\alpha$ such that
$a_\alpha\cap x$ is infinite for each $x\in X$. 
\end{enumerate}

It is straightforward to obtain such a family from,
for example, $\diamondsuit$. $(*)$ gives a strong counterexample
showing that some assumption is needed in Theorem \ref{MA}:

\begin{ex} Assuming $(*)$ there is a compact space $X$ of cardinality $\omega_1$ and of scattered
height 3 such that $X$ is not in the class $\elsig(\loe \omega)$.
\end{ex}

\begin{pf} Let $X$ be the one-point compactification of the
$\Psi$-like space based on the almost disjoint family
$\setof{a_\alpha}{\alpha\in \omega_1}$ witnessing $(*)$, then $X$ is not
in $\elsig(\loe\omega)$. To see this suppose otherwise, and let $C$ be a cover by
second countable (hence countable) compact subsets and $\mathcal
N$ a countable network with respect to $C$. Choose $\beta\in \omega_1$ so that
if $N\in \mathcal N$ and if $N\cap \omega_1$ is countable then $N\cap
\omega_1\subs \beta$. Enumerate the set $\setof{N\cap
(\omega_1\setminus \beta)}{N\in {\mathcal N}}\setminus\{\emptyset\}$
as $\setof{N_k}{k\in \omega}$. 

Choose now a $c\in C$ such that $c\cap (\omega_1\setminus\beta)$ is
not empty. Since $c$ is countable, we may also choose $\alpha$ such
that $a_\alpha\not\in c$ and $a_\alpha\cap N_k$ is infinite for each
$k\in\omega$.  But $a_\alpha \cap c$ is finite, so $U=X\setminus
(\{a_\alpha\}\cup (a_\alpha\setminus c))$ is an open set containing
$c$. However, by choice of $a_\alpha$, there is no element of
${\mathcal N}$ containing $c$ and contained in $U$.
This contradicts that ${\mathcal N}$ is a network at $C$.
\end{pf}


\begin{tw} Assuming $\MAone$, Aronszajn trees are in the class $\elsig(\loe\omega)$.
\end{tw}

\begin{pf} For the topology on $T$, an Aronszajn tree, we define
$[s]=\setof{t\in T}{s\loe t}$ and declare $[s]$ clopen for each $s$ of
successor height in $T$. For each $s\in T$ let $\lceil s
\rceil=\setof{t\in T}{t\loe s}$, the {\it downward closure of
$s$}. The downward closure of each element of $T$ is compact second
countable, so if $\mathcal C=\setof{\lceil s\rceil}{s\in T}$ then
$\mathcal C$ is a cover by second countable compact sets. For $X\subs
T$ we let $\lceil X\rceil=\bigcup\setof{\lceil t \rceil}{t\in X}$.

We will apply $\MAone$ to find a family $\setof{(D_n,F_n)}{n\in \omega}$ such that 

\begin{enumerate}
\item $\lceil D_n\rceil\cap F_n=\emptyset$ for each $n\in\omega$, and
\item for each $t\in T$ and each $F\in [T]^{<\omega}$ such that
$\lceil t \rceil\cap F=\emptyset$, there is $n$ such that $t\in D_n$
and $F\subs F_n$.
\end{enumerate}

Indeed, if such a family exists, it is straightforward to verify that $\setof{\lceil D_n\rceil}{n\in \omega}$ is a network for the cover $\mathcal C$.

To obtain the family, we first define a poset $\mathbb P$ as follows. Let 
$$
{\mathbb P}=\setof{(x,F)}{x,F\in[T]^{<\omega}\mbox{ such that }\lceil x \rceil\cap F=\emptyset}.
$$
We order $\mathbb P$ in the natural way: $(x,F)<(y,G)$ if $y\subs x$
and $F\subs G$. 

\begin{claim} ${\mathbb P}$ is ccc
\end{claim}

\begin{pf} Let $\setof{(x_\alpha,F_\alpha)}{\alpha\in \omega_1}\subs
\mathbb P$. Without loss of generality we may assume that 
\begin{enumerate}
\item[(a)] for $\alpha<\beta$, all elements of $x_\alpha\cup F_\alpha$
lie in levels lower than any element of $x_\beta\cup F_\beta$,
\item[(b)] there is $n,m\in \omega$ such that $|x_\alpha|=n$ and
$|F_\alpha|=m$ for all $\alpha\in \omega_1$. 
\end{enumerate}
We prove by induction on $n$ and $m$ that there is an uncountable
centered subset. The first nontrivial case is where $n=m=1$. In this
case let $x_\alpha=\{s_\alpha\}$ and $F_\alpha=\{t_\alpha\}$. Since we
are assuming $\MAone$ we have that $T$ is special, so we may assume
that 
\begin{enumerate}
\item[(c)] $\setof{t_\alpha}{\alpha\in \omega_1}$ is an antichain in $T$.
\end{enumerate}
If $\setof{\beta}{t_\alpha<s_\beta}$ is countable for each $\alpha$,
then it is straightforward to recursively construct an uncountable
pairwise compatible subset. Otherwise, there is $\alpha_0$ such that
$A=\setof{\beta}{t_{\alpha_0}<s_\beta}$ is uncountable. However, in
this case (a) and (c) imply that $\setof{(s_\beta,t_\beta)}{\beta\in
A}$ is pairwise compatible (in fact centered).

To accomplish the induction step assume that for any family
$\setof{(y_\alpha,G_\alpha)} {\alpha\in \omega_1}$ as above with
$|y_\alpha|=n-1$ has an uncountable centered family and that for any
family $\setof{(y_\alpha,G_\alpha)}{\alpha\in \omega_1}$ with
$|y_\alpha|=n$ and $|G_\alpha|<m$ has an uncountable centered family.

Consider first the case that $m>1$: Fix $t_\alpha\in F_\alpha$ and let
$G_\alpha=F_\alpha\setminus\{t_\alpha\}$. By the inductive assumption
for $m-1$, there is an uncountable $A$ such that
$\setof{(x_\alpha,G_\alpha)}{\alpha\in A}$ is centered. Also by the
induction assumption for $1<m$ there is an uncountable $B\subs A$ such
that $\setof{(x_\alpha,\{t_\alpha\})}{\alpha\in B}$ is centered. By
the definition of the ordering it follows that
$\setof{(x_\alpha,F_\alpha)}{\alpha\in B}$ is centered.

To prove the case where $m=1$, we may assume that $n>1$. Fix
$s_\alpha\in x_\alpha$ and let $y_\alpha=x_\alpha\setminus\{s_\alpha\}$. By
the same procedure as above we may obtain an uncountable $B$ so that both families
$\setof{(y_\alpha,F_\alpha)}{\alpha\in B}$ and
$\setof{(\{s_\alpha\},F_\alpha)}{\alpha\in B}$ are centered. By the
definition of the ordering it follows that
$\setof{(x_\alpha,F_\alpha)}{\alpha\in B}$ is centered.
\end{pf}

To finish the proof of the theorem, take the finite support product of
countably many copies of $\mathbb P$. Take $G$ generic for the family
of dense sets: $D_{t,G}=\setof{p}{(\exists\;n)\; p(n)=(x,F),\; t\in x,
G\subs F}$ where $(t,G)$ range over all pairs $t\in T$ and $G\in
[T]^{<\omega}$ such that $\lceil t\rceil \cap G=\emptyset$. 
Letting 
$$
D_n=\bigcup\{x:\exists p\in G\exists F (p(n)=(x,F)\}
$$
and 
$$
F_n=\bigcup\{F:\exists p\in G\exists x (p(n)=(x,F)\}
$$
it is easy to verify that $\{(D_n,F_n):n\in \omega\}$ satisfy 1.\ and 2.\ as required.

\end{pf}

\section{Open questions}

\begin{question}
Does there exist in ZFC a space in $\kelsigO$ without a dense metrizable subspace?
\end{question}

\begin{question}
Is the class of compact spaces in $\elsig(<\omega)$ absolute?
\end{question}

\begin{question} Does there exist a ($\sig$-compact) space $X$ such
that $\tight(X)>\omega$ and $X\in\elsig(n)$ for some $\ntr$?
\end{question}

\begin{question} Assume that $X^\omega\in\elsigfin$. Is it true in ZFC that
$nw(X)\loe\omega$? \end{question}

\begin{question} \label{fques} Assume that $X\in\elsigO$ and $p\:X\to Y$ is a finite-valued usc mapping, $Y=p(X)$. Must $Y$ belong to $\elsigO$?
\end{question}



\begin{question} Are all Rosenthal compacta in $\elsigO$?
\end{question}

\begin{question} Assume $MA_{\omega_1}$. Is it true that every
scattered compact space of cardinality $\omega_1$ and height $n$,
$n\in\omega$, belongs to $\elsig(\le n+1)$?
\end{question}

\begin{question}\label{qteq} Suppose that $X\in \elsigO$ and $Y$ is a space such that
$C_p(X)$ is homeomorphic to $C_p(Y)$. Must $Y$ belong to $\elsigO$?

\medskip

\end{question}

{\sl Remark.} If $X\in \elsig(<\omega)$ and $C_p(X)$ is homeomorphic
to $C_p(Y)$, then $Y\in \elsig(<\omega)$; this follows from the
invariance of $\elsig(<\omega)$ with respect to finite-valued usc
images and the main theorem in \cite{okt}. By a similar argument, a
positive answer to Question \ref{fques} implies ``yes'' for
Question \ref{qteq}.

\end{document}